# Rate of decay for the mass ratio of pseudo-holomorphic integral 2-cycles

*Costante Bellettini*[*]

Princeton University and IAS

**Abstract**: *We consider any pseudo holomorphic integral cycle in an arbitrary almost complex manifold and perform a blow up analysis at an arbitrary point. Building upon a pseudo algebraic blow up (previously introduced by the author) we prove a geometric rate of decay for the mass ratio towards the limiting density, with an explicit exponent of decay expressed in terms of the density of the current at the point. With a non-explicit exponent this result was proved using different techniques by Pumberger and Rivière in [11].*
  *Primary Subject*: 53C38
  *Secondary Subjects*: 49Q20, 35J99

## 1 Introduction

The possibility that non-smooth behaviour arises for solutions of variational problems or PDEs is a fruitful source of research, both for the understanding of a suitable "weak formulation" of the problem and for the "description of the singularities" of the solutions. We can think for example of harmonic maps, area-minimizing surfaces, mean-curvature flow, hyperbolic conservation laws. A classical tool used to investigate the behaviour at a singular point is the *blow up analysis*, which we now briefly describe in the case of a calibrated cycle, which is relevant for this work.

The *blow up limit* of a current $C$ at a (singular) point $x_0$ requires the following procedure, whose idea goes back to E. De Giorgi [6]. Look at the restriction of $C$ to the ball $B_r(x_0)$ for any small $r$ and then dilate around $x_0$ to the size of the unit ball. When $C$ is a stationary $m$-cycle (in particular when it is calibrated) we have a monotonicity formula whose main consequences are (i) the *mass ratio* $\dfrac{M(C \llcorner B_r(x_0))}{r^m}$ (weakly) decreases as $r \downarrow 0$ (ii) there exist weak limits (as $r \to 0$) of the dilated currents (iii) any such limit is a cone and it is called a *tangent cone* to $C$ at $x_0$. It is an open and difficult problem whether the tangent cone obtained in the limit is *unique* or not: the issue is that by dilating $C$ for different sequences $r_i \downarrow 0$ we might find different cones.

---
[*]*Address*: Centre for Mathematical Sciences, Wilberforce Road, Cambridge CB3 0WB
*E-mail address*: cbellett@dpmms.cam.ac.uk

It is well known (see [7] 5.4.3) that a sufficient condition for the tangent cone to be unique at a point is that the mass ratio at $x_0$ converges "fast enough" to its limit $L(x_0)$ as $r \downarrow 0$. The speed of convergence is often called *rate of decay* for the mass ratio. If there exists a positive increasing function $f$ such that $\int_0^1 \frac{\sqrt{f(r)}}{r} dr < +\infty$ and such that

$$\frac{M(C \llcorner B_r(x_0))}{r^m} - L(x_0) \leq f(r) \tag{1}$$

then the tangent cone at $x_0$ is unique. In particular this will be true if we can prove a rate of decay with $f(r) = Kr^\gamma$ for some $K > 0$ and $\gamma \in (0, 1)$. In this latter case we will say that we have a geometric rate of decay. Let us quote a few proofs where the main goal is to prove the uniqueness of tangents: as we will see in quite a few of them the understanding of the rate of decay is the main tool through which the goal is achieved. In all of the following cases we are not doing justice to the proofs, which contain many interesting and original ideas.

In [17] B. White obtains the uniqueness of tangent cones for mass minimizing two-dimensional cycles. With an over-simplification, we can summarize the steps of his proof as follows: obtain a *Fourier decomposition* of the 1-dimensional cycles defined by intersecting the original current with smaller and smaller spheres centered at $x_0$, get an *epiperimetric inequality* for the current restricted to any small ball around $x_0$, obtain a geometric *rate of decay* for the mass ratio.

In [11] D. Pumberger and T. Rivière prove the uniqueness of tangent cones for semi-calibrated two-dimensional cycles. Roughly speaking again, they argue along the following main steps. Approximate in a suitable sense any semi-calibrated 2-cycle with a pseudo holomorphic one, so that this latter case becomes the main one on which they focus. By means of a *decomposition* into closed Lipschitz curves for the 1-dimensional cycle defined by intersecting the original current with a small sphere centered at $x_0$, they get a Poincaré inequality on this 1-dimensional cycle, from which they obtain a geometric *rate of decay* for the mass ratio.

In other works on (semi-)calibrated 2-cycles proofs of the uniqueness of tangent cones have been given by making a strong use of a property of positive intersection, see the case of integral pseudo-holomorphic 2-cycles in dimension 4 (C. H. Taubes in [16], T. Rivière and G. Tian in [12]) and integral Special Legendrian 2-cycles in dimension 5 (the author and T. Rivière in [3]). This positiveness property is very peculiar to these cases. Moreover on one hand it allows a rather quick proof of the uniqueness, on the other it does not seem easy to obtain a rate of decay for the mass ratio by directly exploiting it.

In [13] the uniqueness for pseudo holomorphic integral 2-dimensional cycles (with a closed two-form locally taming the almost complex structure) is achieved in arbitrary codimension by means of a refinement of the lower-epiperimetric inequality in [17] and by finding a suitable control on the rate of convergence of the mass ratio to its limit.

In [14] L. Simon proves a Lojaciewicz inequality that leads to a rate of decay: in particular it allows to prove that if a tangent cone to a minimal integral current has multiplicity one and has an isolated singularity, then it is unique.

In [1] and [2] the author introduced a new idea to prove uniqueness of tangent cones for pseudo holomorphic cycles of arbitrary dimension. The key step requires to construct a local *pseudo algebraic blow up* of the ambient manifold and of the current at the point (the current obtained after this procedure is called proper tranform). This mimics the classical blow up of singularities of algebraic curves, but must be done by carefully respecting the almost complex structure in a suitable sense, in order not to run into big analysis troubles. Once this step is completed, the proof of the uniqueness is done by working on the proper transform of the original current by arguments that are quite geometric and never require the understanding of the rate of decay. It is however natural to ask whether the rate of decay can be obtained a posteriori and how the pseudo algebraic blow up can simplify the analysis estimates that are needed. In this work we answer these questions by focusing on the case of a two-dimensional pseudo holomorphic integral cycle. As we shall see (Lemma 2.2 and Remark 2.4), the pseudo algebraic blow up on one hand yields an interesting geometric meaning for the left hand side of (1); on the other hand many analysis difficulties are absorbed in the pseudo algebraic blow up and the estimates that need to be done become simpler, in that we are able to get rid of singular terms.

The present work builds a bridge between the approach in [11] and the one in [2]. In both works the analytic tool through which the rate of decay is obtained is the classical hole-filling technique (De Giorgi, Morrey). Many estimates that we will obtain have analogues in [11] and provide interesting geometric interpretations of the arguments there. In particular the Poincaré inequality of [11] turns out to be be a close relative of our Lemma 3.4, with the difference that we will always be working on the proper transform rather than on the original current (see also Remark 3.3 at this stage). We stress that the Poincaré inequality is to some extent the key step in the approach of [11] and it is a very interesting problem to understand if and how this inequality can be generalized to other situations, possibly in higher dimensions (compare the final comment of [11]). As we quickly mentioned earlier, the proof of the Poincaré inequality in [11] relies on the possibility to write a 1-dimensional integral cycle as a sum of closed Lipschitz curves. This decom-position theorem (see [7]) is not available for cycles of dimension higher than 1. In the same vein remark also that [17] requires a Fourier decomposition of 1-dimensional integral cycles and this is again something for which the dimension 1 plays a fundamental role. In the present work we will obtain the analogue of the Poincaré inequality without making use of any decomposition theorem that needs the dimension to be 1. The only decomposition result we will use is available for integral cycles of arbitrary dimension and therefore we get more chances to extend some of the arguments to higher dimensional situations. In this respect, in a work in progress we are investigating a rate of decay for pseudo holomorphic integral cycles of arbitrary dimension. The two-dimensional case deserves however particular attention since, not only it allows an easier presentation of the main ideas, but is also generic from the point of view of the local existence (compare the introduction of [13]). Unlike [11], the geometric rate of decay that we obtain comes with an explicit exponent, depending only on the reciprocal of the density at the point. The present work could thus be used as a starting point for a different proof of the regularity result in [13] (i.e. that pseudo holomorphic integral 2-cycles can only have isolated singularities), which might go through without assuming the existence of a closed two-form taming the almost complex structure[1].

**Setting and main statement.** Let $(\mathcal{M}, J)$ be an almost complex manifold, where the dimension of $\mathcal{M}$ is $(2n + 2)$ and $J$ is the almost complex structure on the tangent bundle. We consider an arbitrary integral cycle of dimension 2 with the property that almost all approximate tangents are positively oriented $J$-invariant 2-planes. Recall that the orientation on $\mathcal{M}$ is induced by $J$ (see e.g. [10]). Such currents are called positive-$(1, 1)$ integral cycles. Recall that the cycle condition (absence of boundary) for an arbitrary $m$-dimensional current $C$ means that it holds, for any compactly supported $(m − 1)$-form $\alpha$, $(\partial C)(\alpha) := C(d\alpha) = 0$.

Such a $T$ can also be viewed as a semi-calibrated cycle, as follows: given $(\mathcal{M}, J)$ as above, it is locally always possible to find a non-degenerate differential form $\omega$ of degree 2 compatible with $J$. The compatibility relies in the fact that $g(\cdot, \cdot) := \omega(\cdot, J\cdot)$ defines a Riemannian metric on $\mathcal{M}$. The tensor $h = g - i\omega$ is called a Hermitian metric on $(\mathcal{M}, J)$. If $d\omega = 0$ then we have a symplectic form, but in general closedness cannot be expected in dimension higher than 4: an example was exhibited on $S^6$ in [5]. The triple $(\mathcal{M}, J, g)$ is an almost Hermitian manifold; when the associated form $\omega$ is closed, we get an almost Kähler manifold. The word "almost" refers to the fact that $J$ can be non-integrable. The form $\omega$ on $(\mathcal{M}, J)$ has (pointwise) unit comass for the associated metric $g(\cdot, \cdot) := \omega(\cdot, J\cdot)$, i.e. it is a semi-calibration. The calibrated 2-planes turn out to be exactly those that are positively oriented $J$-invariant 2-planes, therefore $T$ is a positive-$(1,1)$ integral cycle if and only if it is semicalibrated by $\omega$. It is important to remark that, given an almost complex manifold $(\mathcal{M}, J)$ we can associate a couple $(\omega, g)$ in plenty of ways.

---

[1] A regularity result of this type would extend immediately to arbitrary semi-calibrated integral 2-cycles. Indeed, in view of [4], any semi-calibrated integral 2-cycle can be locally viewed as a pseudo holomorphic one in an almost Hermitian manifold. For the same reason, Theorem 1.1 of this paper can also be extended to arbitrary semi-calibrated integral 2-cycles upon changing the explicit expression of the constant $\gamma$.

In [2] (see also [1]) the author proved, in particular, that such a positive-$(1,1)$ integral cycle $T$ possesses everywhere a unique tangent cone. The notion of tangent cone to a $m$-current $C$ at a point $x_0$ is defined by the following *blow up limit*. Dilate $C$ around $x_0$ of a factor $r$; in normal coordinates around $x_0$ this amounts to pushing forward $C$ via the map $\dfrac{x - x_0}{r}$:

$$(C_{x_0,r} \llcorner B_1)(\psi) := \left[\left(\frac{x - x_0}{r}\right)_* C\right](\chi_{B_1}\psi) = C\left(\chi_{B_r(x_0)}\left(\frac{x - x_0}{r}\right)^* \psi\right). \qquad (2)$$

The almost-monotonicity formula (see [11] Section 2, [1] Prop. 2.1 or [2] Prop. 1 for the semi-calibrated case) gives that the *mass ratio* $\dfrac{M(C \llcorner B_r(x_0))}{\alpha_m r^m}$ (we denote by $\alpha_m$ the volume of the unit $m$-dimensional ball) is monotonically almost-decreasing as $r \downarrow 0$ and therefore, for $r \leq r_0$ (for a small enough $r_0$), we are dealing with a family of cycles $\{C_{x_0,r} \llcorner B_1\}$ in $B_1$ that are equibounded in mass. Thus by Federer-Fleming's compactness theorem (see e.g. [8] page 141 or [7]) we get that there exist weak limits of $C_{x_0,r}$ as $r \to 0$. Every such limit $C_\infty$ is an integer multiplicity rectifiable boundaryless current which turns out to be a cone (recall that a current is said to be a cone with vertex $p$ if it is invariant under homotheties centered at $p$) calibrated by $\omega_{x_0}$ and is called a tangent cone to $C$ at $x_0$. The density of each tangent cone at the vertex is the same as the density of $C$ at $x_0$ (see [9]). It is an open and difficult problem whether the tangent cone obtained in the limit is independent of the sequence of radii yielded by the compactness theorem, i.e. whether the tangent cone at an arbitrary point is *unique* or not.

In the present work, building upon the proof in [2], we consider a positive-$(1,1)$ integral cycle $T$ and we obtain that the mass ratio $\dfrac{M(T \llcorner B_r(x_0))}{\pi r^2}$ decays towards the density $\nu(x_0)$ with a geometric rate. Remark that for semi-calibrated cycles the density is always $\geq 1$ (see the introduction of [2]). In the particular case of dimension 2 we get in addition that the density takes integer values (but we will not need this latter piece of information).

**Theorem 1.1.** *Let $(\mathcal{M}, J)$ be an almost complex manifold and consider in it an arbitrary pseudo holomorphic integral cycle $T$ of dimension $2$. Let $x_0$ be an arbitrary point in the support of $T$. Endow a neighbourhood of $x_0$ with a Riemannian metric compatible[2] with $J$. Denote by $\nu(x_0)$ the density of $T$*

---

[2]As remarked earlier, locally there are infinitely many choices for such a metric. Our result is qualitatively independent of the choice. The actual value of $\tilde{r} > 0$ will depend on this choice.

*at this point. For any $\varepsilon > 0$ there exists $\tilde{r} > 0$ such that for any $r < \tilde{r}$ the following inequality holds:*

$$\frac{M(T \mathbin{\vrule height 1.6ex depth 0pt width 0.4pt \vrule height 0.4pt depth 0pt width 1.2ex } B_r(x_0))}{\pi r^2} - \nu(x_0) \leq K r^\gamma, \tag{3}$$

*with $K$ a small positive constant and $\gamma = \frac{1}{C(1+\varepsilon)\nu(x_0)}$, where $C$ is an explicit universal constant* [3].

**Aknowledgments**. The author wishes to thank T. Colding and C. De Lellis for suggesting the question and T. Rivière and G. Tian for fruitful conversations while the work was in progress.

## 2  Preliminaries.

Theorem 1.1 is a local result, i.e. it only depends on $T$ restricted to a local chart around $x_0$. We therefore assume that we are in the unit ball of $\mathbb{C}^{n+1}$ with a given almost complex structure $J$ and that the point at which we perform the blow up analysis (i.e. the dilations) is the origin itself. We also fix a compatible Riemannian metric $g$ in this chart. By a suitable choice of coordinates we can assume that $J$ coincides with the standard complex structure $J_0$ at the origin and that $g$ coincides with the Euclidean metric at the origin. Define (uniquely) the associated semi-calibration $\omega$ by $g(\cdot, \cdot) = \omega(\cdot, J\cdot)$. The semi-calibration $\omega$ agrees with the standard Kähler form $\omega_0$ at the origin.

For the proof of Theorem 1.1 we will assume (by the result in [2]) that the tangent to $T$ at the origin is unique. In the present work we will build upon the arguments developed in [2] in order to obtain the rate of decay. In this section we set up notations and briefly recall the purpose and effect of the *pseudo-algebraic blow up*, a construction introduced in Sections 3 and 4 of [1] (see also "How to blow up the origin" in Section 3 of [2]). The phrase *pseudo-algebraic blow up* will be kept in order to distinguish this tool from the *blow up analysis*, which is the process of dilations described in (2).

We know that the tangent to $T$ at the origin is a cone of dimension 2 that is moreover $(1,1)$-positive for the standard complex structure $J_0$, therefore it is a finite sum of $J_0$-holomorphic 2-planes, each one counted with an integer multiplicity (for this conclusion the reader may consult e.g. [17], Section 3 of [11] or Section II of [13]): these multiplicities add up to $\nu(x_0)$ (this yields in particular that $\nu(x_0)$ is everywhere integer on $T$).

By standard arguments (e.g. [13] Lemma III.1) that make use of the uniqueness of the tangent cone and the almost monotonicity formula we obtain:

---

[3]Compare (39) for sharper values of $K$ and $\gamma$.

**Lemma 2.1.** *Let $T$ be a positive-$(1,1)$ cycle in $B_1(0)$ and let $D := \oplus N_i [\![D_i]\!]$ be the unique tangent cone at the origin, where the $N_i$'s are positive integers and the $D_i$'s are $J_0$-holomorphic disks through the origin. For any $\varepsilon > 0$ there exists $R > 0$ such that the set*

$$\{0\} \cup \{x \in B_1(0) : dist(x, D) < \varepsilon |x|\}$$

*contains the support of the currents $T_{0,r}$ for any $r < R$.*

In other words, by dilating $T$ enough about the origin, the dilated current lives inside any fixed conical neighbourhood of the tangent cone. This lemma gives the geometric meaning of the tangent cone and of its uniqueness. Moreover remark that the result in Theorem 1.1 is asymptotical, in the sense that we can study the dilated current $T_{0,R}$ (for some $R > 0$) rather than the given $T$; the result obtained for $T_{0,R}$ implies that for $T$. Since $T_{0,R}$ is a positive-$(1,1)$ cycle for the almost complex structure $J_R$ (see "Rescale the foliation" in Section 3 of [1]) obtained by dilating the almost complex structure $J$ in the ball $B_R(0)$ to the unit ball, we can assume in addition that that $J - J_0$ is small in the $C^2$-norm of $B_1(0)$.

In view of this lemma, we can observe that, for each fixed $i$, if $\varepsilon$ is small enough then the current

$$T_{0,R} \llcorner \{x \in B_1(0) : \mathrm{dist}(x, D_i) < \varepsilon |x|\} \tag{4}$$

is still a positive-$(1,1)$ cycle in $B_1(0)$, with respect to the almost complex structure $J_R$ obtained by dilating $J$ of a factor $R$: the tangent cone to the current (4) is clearly the disk $D_i$ counted with multiplicity $N_i$. In this way we can decompose $T_R$ as a sum of positive-$(1,1)$ cycles (that only have the origin as common point of their supports) and we can study these cycles separately and prove Theorem 1.1 for each of them. Then, since $\sum_i N_i = \nu(0)$, it is enough to add up the conclusions (3) obtained for each $i$ to get the theorem for the whole current $T_R$.

*Remark* 2.1. In other words it is enough to prove Theorem 1.1 only in the case that the tangent cone at $x_0$ is a single disk counted with an integer multiplicity and under the further assumption that $T$ is supported in a small conical neighbourhood of this disk (a neighbourhood of the type (4) with $\varepsilon$ as small as we wish).

**The *pseudo-algebraic blow up*.** Let $D$ be the disk supporting the tangent cone to $T$ at the origin. Without loss of generality, let $D$ be the disk $\{z_1 = ... = z_n = 0\}$. By parametrizing the $J_0$-holomorphic disks through the origin with points of $\mathbb{CP}^n$, the disk $D$ corresponds to the point $[0, ..., 0, 1]$ in the homogeneous coordinates of $\mathbb{CP}^n$, or equivalently the origin $(0, .., 0)$ in the chart $z_{n+1} \neq 0$.

By the fixed point techniques described in [12] Lemma A.2, we obtain a pseudo holomorphic polar foliation of a sector that contains $D$. More precisely, we find a family $D^X$ of embedded pseudo holomorphic disks through the origin, whose tangent planes at 0 are parametrized by $X \in \mathbb{CP}^n$ with $X$ ranging in an open ball of $\mathbb{CP}^n$ centered at $[0,...,0,1]$, precisely the ball

$$\mathcal{V} := \left\{ \sum_{j=1}^{n} \frac{|z_j|^2}{|z_{n+1}|^2} < 1 \right\} \subset \mathbb{CP}^n.$$

The family $\{D^X\}$ thus foliates an open set that is a small perturbation of

$$\mathcal{S} = \{(z_1,...z_{n+1}) \in B_1^{2n+2} \subset \mathbb{C}^{n+1} : |(z_1,...,z_n)| < |z_{n+1}|\}.$$

More precisely $\{D^X\}$ is contained in $\{(z_1,...z_{n+1}) \in B_{1+\epsilon}^{2n+2} : |(z_1,...,z_n)| < (1+\epsilon)|z_{n+1}|\}$ and $\{D^X\}$ contains $\{(z_1,...z_{n+1}) \in B_{1-\epsilon}^{2n+2} : |(z_1,...,z_n)| < \frac{1}{1+\epsilon}|z_{n+1}|\}$ for some small $\epsilon > 0$. Remark that $\mathcal{S}$ is the sector foliated by $J_0$-holomorphic disks through the origin, the perturbation is due to the fact that we are taking pseudo holomorphic disks with respect to the almost complex structure $J$, which coincides with $J_0$ at the origin and is a small $C^2$-perturbation of $J_0$ in the unit ball.

By Remark 2.1 the current $T$ is all contained in the foliated sector $\cup_X D^X$, this is why the construction we are performing only needs to be done on the sector rather than on the whole ball $B_1(0)$. Consider the tautological complex line bundle $\widetilde{\mathbb{C}}^{n+1}$ over $\mathbb{CP}^n$. The map $\tilde{\Phi}^{-1}$ constructed in [2] is a diffeomorphism from $\cup_X D^X \setminus \{0\} \approx \mathcal{S}$ into $\widetilde{\mathbb{C}}^{n+1}$. The map $\tilde{\Phi}^{-1}$ is constructed in such a way that the image of the $J$-holomorphic punctured disk $D^X \setminus \{0\}$ is the punctured unit disk in the fiber above $X$, see Figure 1. Denoting by $\mathcal{V} \times \{0\}$ the zero-section of the tautological line bundle $\widetilde{\mathbb{C}}^{n+1}$ restricted to $\mathcal{V}$, we remark that $\left(\tilde{\Phi}^{-1}(\cup_X D^X) \setminus \{0\}\right) \cup (\mathcal{V} \times \{0\})$ is the part of the line bundle lying above $\mathcal{V}$ and we only cover the unit disk in each fiber; we will denote this unit-disk bundle over $\mathcal{V}$ by $\mathcal{A}$. With an abuse of language we can say that the origin is sent, via $\tilde{\Phi}^{-1}$, to the set $\mathcal{V} \times \{0\}$ (this feature justifies the name "blow up" in the sense of algebraic geometry).

By pulling back the almost complex structure $J$ via the projection map $\tilde{\Phi}|_{\mathcal{A} \setminus (\mathcal{V} \times \{0\})} : \mathcal{A} \setminus (\mathcal{V} \times \{0\}) \to (\cup_X D^X) \setminus \{0\}$ we obtain a smooth almost complex structure $I$ on $\mathcal{A} \setminus (\mathcal{V} \times \{0\})$. In Lemma 4.1 of [1] it is proven that $I$ can be defined on the whole of $\mathcal{A}$ by setting it to be the standard almost complex structure of $\widetilde{\mathbb{C}}^{n+1}$ on $\mathcal{V} \times \{0\}$. By making this extension across the zero section $\mathcal{V} \times \{0\}$ we lose the smoothness, indeed $I$ is Lipschitz continuous but not necessarily $C^\infty$. It is moreover a small perturbation of the standard complex structure $I_0$ of $\widetilde{\mathbb{C}}^{n+1}$. In [1], Section 4, we also construct a metric $g$ and a semi-calibration $\theta$ on $\mathcal{A}$ that are Lipschitz continuous on $\mathcal{A}$ and actually smooth away from $\mathcal{V} \times \{0\}$. Together with the almost complex structure $I$ they form a compatible triple $(I, g, \theta)$ in the usual sense, i.e.

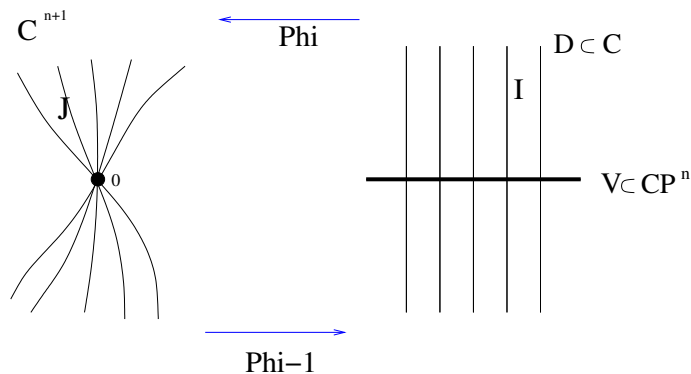

Figure 1: Pseudo-algebraic blow up of the origin.

$g(\cdot, \cdot) = \theta(\cdot, I\cdot)$. Both $g$ and $\theta$ are small perturbations (again in the sense of [1], Lemma 4.1) respectively of the standard metric $g_0$ and the standard Kähler form $\theta_0$ on $\widetilde{\mathbb{C}}^{n+1}$ and actually coincide with $g_0$ and $\theta_0$ on $\mathcal{V} \times \{0\}$.

*Remark* 2.2. Recall that the standard metric on $\widetilde{\mathbb{C}}^{n+1} \subset \mathbb{CP}^n \times \mathbb{C}^{n+1}$ is obtained starting from the Fubini-Study metric on $\mathbb{CP}^n$ and the flat metric on $\mathbb{C}^{n+1}$. Remark in particular that $g$ **is not** the pull back metric via $\tilde{\Phi}$ of the metric that we had in $\mathcal{S}$; this is a key idea in the proof developed in [1] and is connected to the aforementioned fact (see the introduction) that for a given almost complex structure we have plenty of choices for a metric and a non-degenerate two-form so that the compatibility conditions hold.

What we described so far is how to blow up the ambient manifold at the origin in this almost complex setting, by perturbing the classical blow up of $\mathbb{C}^{n+1}$ used in algebraic and symplectic geometry. The purpose of this construction is to push-forward the current $T$ (that is supported inside $\mathcal{S}$ by assumption) via $\tilde{\Phi}^{-1}$ and obtain a new positive-$(1,1)$ cycle $P$ in $\mathcal{A}$. Let us sketch how to do it. The new cycle $P$ will be called *proper transform* of $T$.

We can surely perform a "partial" push-forward, i.e. we remove a small neighbourhood of the origin and push-forward the current $T \llcorner \left( \mathcal{S} \setminus \overline{B_r(0)} \right)$. The map $\tilde{\Phi}^{-1}$ is well-defined, smooth and pseudo holomorphic with respect to $J$ and $I$. In order to define the push-forward on the whole sector we would like to send $r \to 0$ and take the limit of the partial push-forwards. The hard issue to deal with when $r$ goes to 0 is that the map $\tilde{\Phi}^{-1}$ degenerates at the origin. First of all it is not clear that we can obtain $P$ in this way, since in order to pass to the limit we need a uniform bound (indipendently of $r$) on the masses of the push-forwards of $T \llcorner \left( \mathcal{S} \setminus \overline{B_r(0)} \right)$. The second difficulty is that by cutting out $B_r(0)$ we are creating boundary and it is not clear whether this boundaries disappear after having taken the limit of the partial push-forwards: it is a prori conceivable that we might find some boundary on $\mathcal{V} \times \{0\}$ in the limit. These issues are dealt with in [1] Lemmas 4.2 and

4.3 (see also [2] Lemma 3.3), where we show that, by sending $r$ to 0, the push-forward
$$P := \lim_{r \to 0} \left(\tilde{\Phi}^{-1}\right)_* \left[T \llcorner \left(\mathcal{S} \setminus \overline{B_r(0)}\right)\right]$$
is well defined as a boundaryless integral current in $\mathcal{A}$ of finite mass. The fact that $P$ is positive-$(1,1)$ with respect to $I$ follows quite easily, since the map we are using is pseudo holomorphic. The proper tranform $P$ is thus semi-calibrated by $\theta$ in the manifold $(\mathcal{A}, g)$.

*Remark* 2.3. Using the splitting $\mathbb{C}^n \times \mathbb{C}$, keeping in mind that the "vertical fiber $\mathbb{C}$ is a $I$-pseudo holomorphic 2-plane[4], we can see that the almost complex structure $I$ has the form

$$\begin{pmatrix} I_1 & 0 \\ E & I_2 \end{pmatrix}$$

with $I_1^2 = -Id_{2n}$, $I_2^2 = -Id_2$ and $E$ small in $C^0$-norm. More precisely the latter norm it is controlled by $|I - I_0|(p) \leq \varepsilon \, \text{dist}(p, \mathcal{V} \times \{0\})$ (see Lemma 4.1 of [1]). For any vector $\vec{\tau}$, with respect to the splitting $\mathbb{C}^n \times \mathbb{C}$ for the tangent space at the base-point of $\vec{\tau}$, we decompose $\vec{\tau} = \vec{\tau}_h + \vec{\tau}_v$. We have, denoting by $|\cdot|$ the standard flat metric,

$$\langle d\eta_j, \vec{\tau} \wedge I_0 \vec{\tau} \rangle - \varepsilon \, |\vec{\tau}_h|^2 \leq \langle d\eta_j, \vec{\tau} \wedge I \vec{\tau} \rangle \leq \langle d\eta_j, \vec{\tau} \wedge I_0 \vec{\tau} \rangle + \varepsilon \, |\vec{\tau}_h|^2. \quad (5)$$

Since $d\eta$ is the area form for the flat metric, we know that

$$\left\langle \sum_{j=1}^n d\eta_j, \vec{\tau} \wedge I_0 \vec{\tau} \right\rangle = \langle d\eta, \vec{\tau} \wedge I_0 \vec{\tau} \rangle = |\vec{\tau}_h|^2, \quad (6)$$

so that (5) becomes

$$(1 - n\,\varepsilon) \sum_{j=1}^n \langle d\eta, \vec{\tau} \wedge I_0 \vec{\tau} \rangle \leq \langle d\eta, \vec{\tau} \wedge I \vec{\tau} \rangle \leq (1 + n\,\varepsilon) \langle d\eta, \vec{\tau} \wedge I_0 \vec{\tau} \rangle. \quad (7)$$

In particular $\langle d\eta, \vec{\tau} \wedge I\vec{\tau} \rangle$ is always non-negative as long as $\|I - I_0\|_\infty$ is small enough. This is true as long as $T$ is dilated enough about the origin before the procedure of the pseudo holomorphic blow up begins (see Remark 2.1).

---

[4]Without this assumption the estimate we would get instead of (5) would be $\langle d\eta_j, \vec{\tau} \wedge I_0\vec{\tau} \rangle - \varepsilon \, |\vec{\tau}|^2 \leq \langle d\eta_j, \vec{\tau} \wedge I\vec{\tau} \rangle \leq \langle d\eta_j, \vec{\tau} \wedge I_0\vec{\tau} \rangle + \varepsilon \, |\vec{\tau}|^2$, not sufficient for our purposes in view of (6).

**A geometric meaning for the left-hand side of (3)**. Now we will see how a geometric idea as the *pseudo algebraic blow up* that we outlined can have a direct impact on the analysis involved in getting the rate of decay of Theorem 1.1.

We need a digression on the construction of $\tilde{\Phi}$ done in [1], Sections 3 and 4. The map $\tilde{\Phi}$ is obtained via a composition $\Psi \circ \Phi$, where $\Phi$ is the standard projection map from the line bundle $\widetilde{\mathbb{C}}^{n+1}$ to $\mathbb{C}^{n+1}$ and $\Psi$ is a $C^2$ diffeomorphism whose key property relies in the fact that it fixes the origin and it sends flat $J_0$-holomorphic disks through the origin of $\mathbb{C}^n$ to embedded $J$-pseudo holomorphic disks through the origin and in doing so it provides a polar foliation of $\mathcal{S}$. The extra properties of $\Psi$ is that $|D\Psi - \mathbb{I}|(p) \leq Cr$, where $r$ is the distance of $p$ to the origin, $\mathbb{I}$ is the identity matrix and $C = \|D^2\Psi\|_\infty$. So we can think of $\Psi$'s differential at the origin as being $\mathbb{I}$. Just like in [1], let us think of the current $\left(\Psi^{-1}\right)_* T$: it is semicalibrated by the two-form $\Psi^*\omega$ with respect to the metric $\Psi^*g$. Denote by $\tilde{B}_r(0)$ the image via $\Psi$ of the (Euclidean) ball $B_r(0)$ and choose $r^-$ and $r^+$ so that $\overline{B_{r^-}(0)}$ and $B_{r^+}(0)$ are the best possible approximations of $B_r(0)$ respectively from inside and from outside. Then

$$\frac{r^{-2}}{r^2} \frac{M(T \llcorner B_{r^-}(0))}{r^{-2}} \leq \frac{M(T \llcorner \tilde{B}_r(0))}{r^2} \leq \frac{r^{+2}}{r^2} \frac{M(T \llcorner B_{r^+}(0))}{r^{+2}}. \quad (8)$$

The property $|D\Psi - \mathbb{I}|(p) \leq Cr$ then implies that the inequalities (8) become

$$(1-Cr)^2 \frac{M(T \llcorner B_{r^-}(0))}{r^{-2}} \leq \frac{M(T \llcorner \tilde{B}_r(0))}{r^2} \leq (1+Cr)^2 \frac{M(T \llcorner B_{r^+}(0))}{r^{+2}}. \quad (9)$$

The latter inequalities mean that we can study the mass ratio $\frac{M(T \llcorner \tilde{B}_r(0))}{r^2}$ and prove the rate of decay for this expression rather than for the original one. In other words we can study the ratio

$$\frac{\left((\Psi_*^{-1}T) \llcorner B_r(0)\right)(\Psi^*\omega)}{\pi r^2},$$

where $B_r$ is an Euclidean ball.

We need a further observation that allows a nice approximation of the left hand side of (3). For the following discussion denote by $T^\Psi$ the current $\Psi_*^{-1}T$. From the proof of Theorem 5.6 in [9] we have

$$\frac{\left(T^\Psi \llcorner B_r(0)\right)(\omega_0)}{\pi r^2} - \frac{\left(T^\Psi \llcorner B_s(0)\right)(\omega_0)}{\pi s^2} = \int_{B_r(0) \setminus B_s(0)} \frac{\langle (\omega_0)_t, \vec{T}^\Psi \rangle}{|x|^2} d\|T^\Psi\|, \quad (10)$$

where $\|T^\Psi\|$ is the Radon measure associated to the current $T^\Psi$, $\vec{T}^\Psi$ stands for the oriented approximate tangent plane to $T^\Psi$ represented as a unit simple 2-vector ($\|T^\Psi\|$-almost everywhere well-defined) and $(\omega_0)_t$ is the tangential part of the parallel form $\omega_0$ as defined in [9] (5.6), i.e.

$$(\omega_0)_t := \frac{\partial}{\partial r} \lrcorner (dr \wedge \omega_0). \tag{11}$$

Remark that, since $T^\Psi$ is semi-calibrated by $\Psi^*\omega$ (and $\Psi^*\omega = \omega_0$ at the origin), we have $M(T \llcorner B_r(0)) = (T \llcorner B_r(0))(\Psi^*\omega) = (T \llcorner B_r(0))(\omega_0) + (T \llcorner B_r(0))(\Psi^*\omega - \omega_0)$ and the comass of $\Psi^*\omega - \omega_0$ in $B_r(0)$ is controlled by $Cr$, where $C = \|\nabla\omega\|_\infty \cdot \|D^2\Psi\|_\infty$. Thus we get

$$(1 - Cr)M\left(T^\Psi \llcorner B_r(0)\right) \leq \left(T^\Psi \llcorner B_r(0)\right)(\omega_0) \leq (1 + Cr)M\left(T^\Psi \llcorner B_r(0)\right), \tag{12}$$

so the ratios $\frac{(T^\Psi \llcorner B_r(0))(\omega_0)}{\pi r^2}$ and $\frac{M(T^\Psi \llcorner B_r(0))}{\pi r^2}$ have the same limits as $r \downarrow 0$. Send now $s \downarrow 0$ in (10) and use (12) for the first term, we obtain

$$(1 - Cr)\frac{M\left(T^\Psi \llcorner B_r(0)\right)}{\pi r^2} - \nu(0) \leq \int_{B_r(0)} \frac{\langle (\omega_0)_t, \vec{T}^\Psi \rangle}{|x|^2} d\|T^\Psi\|, \tag{13}$$

and thus, upon taking $r$ small enough to ensure $\frac{M(T^\Psi \llcorner B_r(0))}{\pi r^2} \leq 2\nu(0)$,

$$\frac{M\left(T^\Psi \llcorner B_r(0)\right)}{\pi r^2} - \nu(0) \leq \int_{B_r(0)} \frac{\langle (\omega_0)_t, \vec{T}^\Psi \rangle}{|x|^2} d\|T^\Psi\| + 2\nu(0)Cr. \tag{14}$$

For the proof of Theorem 1.1 we therefore need to control the term $\int_{B_r(0)} \frac{\langle (\omega_0)_t, \vec{T}^\Psi \rangle}{|x|^2} d\|T^\Psi\|$, since the remaining term $2\nu(0)Cr$ can be absorbed by $Kr^\gamma$ simply by taking $\tilde{r}$ small enough.

For the estimate of $\int_{B_r(0)} \frac{\langle (\omega_0)_t, \vec{T}^\Psi \rangle}{|x|^2} d\|T^\Psi\|$ we are going to interpret this term geometrically, in particular we need to understand how to interpret the action of $T^\Psi$ on the two-form $\frac{(\omega_0)_t}{|x|^2}$ (which has a singularity at the origin).

The pseudo algebraic blow up recalled in Section 2 will be the key tool to this geometric interpretation. Denote in the sequel with $\omega_{\mathbb{CP}^n}$ the standard Kähler form on $\mathbb{CP}^n$ pulled back to $\mathcal{A}$ via the tautological projection from $\mathcal{A}$ to $\mathcal{V}$. We will also use the notation $\mathcal{A}_r$ to indicate the inverse image of $B_r(0)$ via the map $\Phi$. To avoid confusion, remark that the proper transform of $T$ via $\tilde{\Phi} := \Psi \circ \Phi$ is defined as the proper transform of $T^\Psi$ via $\Phi$.

**Lemma 2.2.** *Let $P$ be the proper transform of $T$. Then*

$$\int_{B_r(0)} \frac{\langle (\omega_0)_t, \vec{T}^\Psi \rangle}{|x|^2} d\|T^\Psi\| = (P \llcorner \mathcal{A}_r)(\omega_{\mathbb{CP}^n}). \tag{15}$$

*proof of Lemma 2.2.* The explicit expression of $\Phi$, the standard projection map associated to the standard algebraic blow up of the origin gives with a straightforward computation that $\Phi^*\left(\frac{(\omega_0)_t}{|x|^2}\right) = \omega_{\mathbb{CP}^n}$. □

*Remark* 2.4. The significance of this lemma is evident: while on the left-hand side we have $T$ acting on a form that has a singularity at the origin, on the right-hand side we have the action of the integral 2-cycle $P$ on a smooth form. We got rid of the singularity and this is a big advantage for the analysis estimates that are needed to achieve the rate of decay. Basically the difficulties that are present in estimating the action on the singular form on the left hand side have been absorbed in the process of the *pseudo holomorphic blow up* by showing that $P$ is indeed well-defined as an integral cycle.

## 3 Proof of the rate of decay.

In order to understand the rate of decay for the left-hand side of (15) we will use De Giorgi/Morrey's technique. Some preliminary considerations (that simplify the task) are in order.

Recall that the pseudo algebraic blow up is done in a sector $\mathcal{S}$ that is a conical neighbourhood of the tangent 2-plane to the cycle $T$ and moreover $\mathcal{S}$ contains the support of $T$. Therefore, as we saw in Section 2, the proper transform $P$ is supported in the image of $\mathcal{S}$ via $\Phi^{-1}$, i.e. in the unit-disk bundle over $\mathcal{V}$ that we denoted by $\mathcal{A}$. With this in mind, when we estimate the right hand side of (15) we can actually work in the chart that describes the unit-disk bundle $\mathcal{A}$ over $\mathcal{V}$ as a subset of $\mathbb{C}^n \times \mathbb{C}$. Use coordinates $(\zeta_1, ..., \zeta_n, Z)$, $\zeta_j = x_j + iy_j$, on $\mathbb{C}^n \times \mathbb{C}$. The two-form $\omega_{\mathbb{CP}^n}$ (the standard Kähler form on $\mathbb{CP}^n$ extended to the cartesian product $\mathbb{C}^n \times \mathbb{C}$ independently of the variable $Z$) reads $\frac{\sum_{i=1}^n dx^i \wedge dy^i}{(1 + |\vec{z}|^2)^{1/2}}$ in this chart. We actually know (compare Remark 2.1) that the support of $T$ is contained in a smaller sector of amplitude as small as we wish, say $\varepsilon$, i.e. in $\{(z_1, ...z_{n+1}) \in B_1^{2n+2} \subset \mathbb{C}^{n+1} : |(z_1, ..., z_n)| < \varepsilon |z_{n+1}|\}$: the proper transform of $T$ is thus supported in the image via $\Phi^{-1}$ of this smaller sector, that is an open set of the form $B_\varepsilon^{2n} \times \mathbb{C}$ inside $\mathbb{C}^n \times \mathbb{C}$, where $B_\varepsilon^{2n}$ is an open ball centered at the origin of $\mathbb{C}^n$ with radius $\varepsilon$. This means that when we estimate the right hand side of (15) working in the chart $\mathcal{V} \times \mathbb{C}$ we can actually estimate the action of $P$ on the standard symplectic form $\omega_{\mathbb{C}^n} = \sum_{i=1}^n dx^i \wedge dy^i$ rather than on $\omega_{\mathbb{CP}^n}$.

Indeed by the considerations just made it holds

$$(P \llcorner \mathcal{A}_r)(\omega_{\mathbb{CP}^n}) \leq (1+\epsilon)(P \llcorner \mathcal{A}_r)(\omega_{\mathbb{C}^n}),$$

with $\epsilon$ independent of $r$, so we can easily see that once we will have obtained $(P \llcorner \mathcal{A}_r)(\omega_{\mathbb{C}^n}) \leq Kr^\gamma$ the same decay will follow for $(P \llcorner \mathcal{A}_r)(\omega_{\mathbb{CP}^n})$. We therefore start now the proof of this decay estimate through De Giorgi/Morrey's technique This requires to bound the quantity $(P \llcorner \mathcal{A}_R)(\omega_{\mathbb{C}^n})$ with a constant (independent of $R$) times the quantity $(P \llcorner (\mathcal{A}_{2R} \setminus \mathcal{A}_R))(\omega_{\mathbb{C}^n})$. The result will be achieved at the end of this section.

Consider on $\mathbb{C}^n$ the one-form $\eta := \sum_{j=1}^{n} x_j dy^j - y_j dx^j$. We will also use the notations $\eta_1 := x_1 dy^1 - y_1 dx^1, ..., \eta_n = x_n dy^n - y_n dx^n$ and the complex variables $\zeta_j = x_j + iy_j$. In view of the estimate on $(P \llcorner \mathcal{A}_R)(\omega_{\mathbb{C}^n})$ that we are seeking we observe that since $P$ is a cycle it holds:

$$(P \llcorner \mathcal{A}_R)(\omega_{\mathbb{C}^n}) = \langle P, |Z| = R \rangle(\eta).$$

Motivated by this, in the lemmas to come we will be dealing with estimates for $\langle P, |Z| = R \rangle(\eta)$.

*Remark* 3.1. The integral cycle $P$ is semi-calibrated and as such satisfies an almost monotonicity formula for the mass ratio at *every* point. It is well-known (compare the introduction of [2]) that this yields that the density function can be well-defined everywhere and is upper semi-continuous, in particular it takes values $\geq 1$ everywhere on the support of $P$. In the following we will denote this density function by $f$. These considerations yield in particular that the rectifiable set $\mathcal{P}$ underlying the current $P$ agrees a.e. with its closure, in other words we can choose (and will do so throughout) a representative for $\mathcal{P}$ that is **closed**.

We will be working with slices of $P \llcorner (\mathcal{A}_{2R} \setminus \mathcal{A}_R)$ with respect to the slicing function $|Z| = r$ for $r \in [R, 2R]$. By abuse of language we will say that we are slicing $P$ with tubes of the form $\mathbb{C}^n \times S^1_r$.

Denote by $\vec{t}$ the unit tangent vectors to the "vertical circles" $S^1_r$. The gradient of the slicing function $|Z|$ is the vector $I_0 \vec{t}$, therefore the Jacobian relative to $P$ for this slicing function is $|\mathrm{proj}_{\vec{P}}(I_0 \vec{t})|$, the length of the projection of the gradient onto the approximate tangent $\vec{P}$. We will denote by $\vec{\tau}$ the approximate (unit) tangent to the slice $\langle P, |Z| = r \rangle$ ($\mathcal{H}^1$-a.e. well-defined on each slice). Since $P$ is a pseudo holomorphic cycle the approximate tangent $\vec{P}$ is $\mathcal{H}^2$-a.e. given by $\vec{\tau} \wedge I\vec{\tau}$.

With this in mind, we use the coarea formula (see e.g. [8] Section 2.1.5 Theorem 3)

$$\int_{R}^{2R} \left( \int_{\mathcal{P} \cap \{|Z|=r\}} \frac{(\vec{\tau} \wedge I\vec{\tau})(d\eta)}{|\mathrm{proj}_{\vec{P}}(I_0\vec{t})|} f d\mathcal{H}^1 \right) dr =$$
$$\int_{\mathcal{P} \llcorner (\mathcal{A}_{2R} \setminus \mathcal{A}_R)} (\vec{\tau} \wedge I\vec{\tau})(d\eta) f d\mathcal{H}^2 = (P \llcorner (\mathcal{A}_{2R} \setminus \mathcal{A}_R))(d\eta) \quad (16)$$

on the rectifiable set $\mathcal{P}$ for the non-negative (see Remark 2.3) integrand $f \frac{(\vec{\tau} \wedge I\vec{\tau})(d\eta)}{|\mathrm{proj}_{\vec{P}}(I_0\vec{t})|}$, where $f \geq 1$ is the density function associated to $\mathcal{P}$ for the current $P$, and the slicing function that we are using is $|Z| = r$ for $r \in [R, 2R]$. Since $\mathcal{P}$ is closed by Remark 3.1, the set $\mathcal{P} \cap \{|Z| = r\}$ is also closed and from [8] Section 2.1.5 Proposition 2 we also know that it is countably $\mathcal{H}^1$-rectifiable for a.e. $r$: these two facts together yield that for a.e. $r$ the set $\mathcal{P} \cap \{|Z| = r\}$ is a closed $\mathcal{H}^1$-rectifiable set.

Remark 3.2. Observing that if $(\vec{\tau} \wedge I\vec{\tau})(d\eta) = 0$ then $\mathrm{proj}_{\vec{P}}(I_0\vec{t}) \neq 0$, we can see that equality (16) gives in particular that for a.e. $r$

$$\int_{\mathcal{P} \cap \{|Z|=r\}} \frac{(\vec{\tau} \wedge I\vec{\tau})(d\eta)}{|\mathrm{proj}_{\vec{P}}(I_0\vec{t})|} f d\mathcal{H}^1$$

is finite and thus $|\mathrm{proj}_{\vec{P}}(I_0\vec{t})| = 0$ only on a $\mathcal{H}^1$-null subset of $\mathcal{P} \cap \{|Z| = r\}$.

By slicing theory (see [8] Section 2.2.5) a.e. slice $\langle P, |Z| = r \rangle$ is well-defined as an integral 1-current with underlying $\mathcal{H}^1$-rectifiable set given by

$$\mathcal{P} \cap \{|Z| = r\} \setminus \{|\mathrm{proj}_{\vec{P}}(I_0\vec{t})| = 0\}$$

and multiplicity function given by $f|_{\mathcal{P} \cap \{|Z|=r\} \setminus \{|\mathrm{proj}_{\vec{P}}(I_0\vec{t})|=0\}}$.

By the observation just made we then get that for a.e. $r$ the set

$$\mathcal{P} \cap \{|Z| = r\}$$

is a good (**closed**) representative for the $\mathcal{H}^1$-rectifiable set underlying the slice $\langle P, |Z| = r \rangle$ and $f|_{\mathcal{P} \cap \{|Z|=r\}}$ is the associated ($\mathbb{N}$-valued) multiplicity.

A further strightforward and important consequence of equality (16), which we can now rewrite

$$\int_{R}^{2R} \left( \int_{\langle P, |Z|=r \rangle} \frac{(\vec{\tau} \wedge I\vec{\tau})(d\eta)}{|\mathrm{proj}_{\vec{P}}(I_0\vec{t})|} d\|\langle P, |Z| = r \rangle\| \right) dr =$$
$$\int_{P \llcorner (\mathcal{A}_{2R} \setminus \mathcal{A}_R)} (\vec{\tau} \wedge I\vec{\tau})(d\eta) d\|P\| = (P \llcorner (\mathcal{A}_{2R} \setminus \mathcal{A}_R))(d\eta), \quad (17)$$

is the fact that we can choose "good slices" $\langle P, |Z| = r \rangle$ in the following sense:

**Lemma 3.1.** *For at least half of the slices $\langle P, |Z| = \rho \rangle$ with $\rho \in [R, 2R]$ it holds*

$$\int_{\langle P,|Z|=\rho\rangle} \frac{(\vec{\tau} \wedge I\vec{\tau})\,(d\eta)}{|proj_{\vec{P}}(I_0\vec{t})|} d\|\langle P,|Z|=\rho\rangle\| \leq \frac{2}{R} \left(P \llcorner (\mathcal{A}_{2R} \setminus \mathcal{A}_R)\right)(d\eta). \quad (18)$$

**Decomposition of the slices**. Any slice of the form $\langle P, |Z| = r \rangle$ is an integral 1-cycle and by the decomposition theorem in [7] 4.2.25 we know that $\langle P, |Z| = r \rangle$ can be decomposed in sub-cycles in a very strong sense, i.e. $\langle P, |Z| = r \rangle$ can be represented as a sum of integral 1-cycles each of which can be parametrized as a Lipschitz curve and moreover $M(\langle P, |Z| = r \rangle)$ equals the sum of the masses of these sub-cycles. From this we easily deduce the following weaker decomposition result, which will be used in the sequel. There exist integral 1-cycles $N_i$'s such that (whenever the slice exists)

$$\langle P, |Z| = r \rangle = \sum_{j=1}^{\infty} N_j \text{ with } \sum_{j=1}^{\infty} M(N_j) = M(\langle P, |Z| = r \rangle), \quad (19)$$

where each $N_j$ is an integral 1-cycle supported on a 1-rectifiable set $\mathcal{N}_i$ that is compact and connected and the $\mathcal{N}_i$'s form a disjoint collection.

In other words for a.e. $Z$ we have $\mathcal{P} \cap \{|Z| = r\} = \sqcup_{i=1}^{\infty} \mathcal{N}_i$ in the a.e. sense and $N_i$ is the current of integration on $\mathcal{N}_i$ with multiplicity $f|_{\mathcal{N}_i}$. We prefer to use the weaker decomposition result in view of possible extension of the approach in this paper to higher dimensional situations, where the stronger decomposition is not available. In the appendix we show that the weaker decomposition result can be proven directly using only elementary tools that would be available in any dimension.

**Lemma 3.2.** *Let $N$ be an integral cycle of dimension 1 in $\mathbb{C}^n \times S^1_\rho$, where $S^1_\rho$ stands for the circle of radius $\rho$. Let $j \in \{1,...,n\}$. Assume that the underlying 1-rectifiable set $\mathcal{N}$ to the current $N$ is a compact and connected set and that the image of $\mathcal{N}$ via the projection $\mathbb{C}^n \times S^1_\rho \to \mathbb{C}^n$ contains the origin. Denote by $\vec{\tau}$ the unit vector that is the approximate tangent to $N$ and by $\vec{t}$ the unit tangent vector to the fiber $S^1_\rho$. The 2-plane spanned by the vectors $\vec{\tau}$ and $I\vec{\tau}$ is denoted by $\vec{\tau} \wedge I\vec{\tau}$. Then, for any $j \in \{1,...,n\}$ we have*

$$|N(\eta_j)| \leq \frac{1}{2}|\mu| \int_{\mathcal{N}} \frac{\langle d|\zeta_j|, \vec{\tau}\rangle^2 + \left\langle \frac{\eta_j}{|\zeta_j|}, \vec{\tau}\right\rangle^2}{|proj_{\vec{\tau} \wedge I\vec{\tau}}(I_0\vec{t})|} d\|N\|. \quad (20)$$

*The constant $|\mu|$ on the right-hand side is the total mass of the measure $\mu := |proj_{\vec{P}}(I_0\vec{t})| \, \|N\|$, i.e. the measure associated to the current $N$ weighted with the positive factor $|proj_{\vec{P}}(I_0\vec{t})|$.*

***proof of Lemma 3.2.*** Without loss of generality we will take $j = 1$.

$$|N(\eta_1)| = \left| \int_{\mathcal{N}} \langle \eta_1, \vec{\tau} \rangle \, d\|N\| \right| \leq$$

$$\leq \int_{\mathcal{N}} \frac{\left| \left\langle \frac{\eta_1}{|\zeta_1|}, \vec{\tau} \right\rangle \right|}{|\text{proj}_{\vec{\tau} \wedge I\vec{\tau}}(I_0\vec{t})|} |\text{proj}_{\vec{\tau} \wedge I\vec{\tau}}(I_0\vec{t})|^{1/2} |\zeta_1| |\text{proj}_{\vec{\tau} \wedge I\vec{\tau}}(I_0\vec{t})|^{1/2} d\|N\| \leq \quad (21)$$

$$\leq \left( \int_{\mathcal{N}} \frac{\left| \left\langle \frac{\eta_1}{|\zeta_1|}, \vec{\tau} \right\rangle \right|^2}{|\text{proj}_{\vec{\tau} \wedge I\vec{\tau}}(I_0\vec{t})|} d\|N\| \right)^{\frac{1}{2}} \left( \int_{\mathcal{N}} |\zeta_1|^2 |\text{proj}_{\vec{\tau} \wedge I\vec{\tau}}(I_0\vec{t})| d\|N\| \right)^{\frac{1}{2}},$$

where we used Cauchy-Schwarz inequality in the last step. We proceed to an important estimate for the second factor in the last line. First of all we need the following analog of the fundamental theorem of calculus:

$$\max_{\mathcal{N}} |\zeta_1| \leq \int_{\mathcal{N}} |\langle d|\zeta_1|, \vec{\tau} \rangle| \, d\|N\|. \quad (22)$$

In order to prove (22) we observe that $|\langle d|\zeta_1|, \vec{\tau} \rangle|$ is the Jacobian relative to the rectifiable set supporting $\mathcal{N}$ of the map $|\zeta_1| : \mathcal{N} \to \mathbb{R}$. We therefore write by means of the coarea formula (see [8] Section 2.1.5 ):

$$\int_{\mathcal{N}} |\langle d|\zeta_1|, \vec{\tau} \rangle| \, d\|N\| = \int_0^{\max_{\mathcal{N}} |\zeta_1|} \left( \int_{|\zeta_1|=s} 1 \right) ds \geq \max_{\mathcal{N}} |\zeta_1|, \quad (23)$$

where in the last step we are using that there is at least a point in each set $|\zeta_1|^{-1}(s)$ by the assumption that $\mathcal{N}$ is connected and by the fact that $\mathcal{N}$ was centered with respect to the origin of $\mathbb{C}^n$. Rewrite now (22) as

$$\max_{\mathcal{N}} |\zeta_1| \leq \int_{\mathcal{N}} \frac{|\langle d|\zeta_1|, \vec{\tau} \rangle|}{|\text{proj}_{\vec{\tau} \wedge I\vec{\tau}}(I_0\vec{t})|} |\text{proj}_{\vec{\tau} \wedge I\vec{\tau}}(I_0\vec{t})| d\|N\| = \int_{\mathcal{N}} \frac{|\langle d|\zeta_1|, \vec{\tau} \rangle|}{|\text{proj}_{\vec{\tau} \wedge I\vec{\tau}}(I_0\vec{t})|} d\mu, \quad (24)$$

with the measure $\mu := |\text{proj}_{\vec{\tau} \wedge I\vec{\tau}}(I_0\vec{t})| \, \|N\|$ (as in the statement Lemma 3.2). By Jensen's inequality we obtain

$$\left( \max_{\mathcal{N}} |\zeta_1| \right)^2 \leq |\mu| \int_{\mathcal{N}} \frac{\langle d|\zeta_1|, \vec{\tau} \rangle^2}{|\text{proj}_{\vec{\tau} \wedge J\vec{\tau}}(J_0\vec{t})|^2} d\mu. \quad (25)$$

Now by integrating the latter inequality with respect to the measure $\mu$ we infer the desired bound on the last factor of (21):

$$\int_{\mathcal{N}} |\zeta_1|^2 |\text{proj}_{\vec{\tau} \wedge I\vec{\tau}}(I_0\vec{t})| d\|N\| = \int_{\mathcal{N}} |\zeta_1|^2 d\mu \leq$$

$$\leq |\mu|^2 \int_{\mathcal{N}} \frac{\langle d|\zeta_1|, \vec{\tau} \rangle^2}{|\text{proj}_{\vec{\tau} \wedge I\vec{\tau}}(I_0\vec{t})|^2} d\mu = |\mu|^2 \int_{\mathcal{N}} \frac{\langle d|\zeta_1|, \vec{\tau} \rangle^2}{|\text{proj}_{\vec{\tau} \wedge I\vec{\tau}}(I_0\vec{t})|} d\|N\|. \quad (26)$$

Combining (21) with (26) and then using the inequality $\sqrt{ab} \leq \frac{1}{2}(a+b)$ (valid for any $a, b \geq 0$) we conclude the proof of the Lemma. □

Next we are going to relate $\langle d|\zeta_j|, \vec{\tau}\rangle^2 + \left\langle \frac{\eta_j}{|\zeta_j|}, \vec{\tau} \right\rangle^2$, which appears in the right-hand side of (20) to the expression $(\vec{\tau} \wedge I\vec{\tau})(d\eta_j)$. The significance of this will appear when we will apply lemma 3.2 to a 1-dimensional slice of $P$, namely the slice of $P$ with a tube of the form $\mathbb{C}^n \times S^1_\rho$: indeed in that case $\vec{\tau}$ is a unit vector in the approximate tangent of $P$, therefore, since $P$ is pseudo holomorphic with respect to $I$, $\vec{\tau} \wedge I\vec{\tau}$ represents the approximate tangent itself.

**Lemma 3.3.** *Let $\vec{\tau}$ be a unit vector in $\mathbb{C}^n \times S^1_\rho \subset \mathbb{C}^n \times \mathbb{C}$ and let $I$ be an almost complex structure on $\mathbb{C}^n \times \mathbb{C}$ with the properties that $I$ is the standard complex structure on $\mathbb{C}^n \times \{0\}$ and each 2-plane $\{(\zeta_1, ..., \zeta_n)\} \times \mathbb{C}$ is pseudo-holomorphic. Let $\|I - I_0\|_{C^0} \leq \frac{\varepsilon}{n}$, for a certain $\varepsilon \geq 0$ (we will use this lemma when the $C^0$-norm of $I - I_0$ is small). Then we have*

$$(1-\varepsilon) \sum_{j=1}^n \left( \langle d|\zeta_j|, \vec{\tau}\rangle^2 + \left\langle \frac{\eta_j}{|\zeta_j|}, \vec{\tau} \right\rangle^2 \right) \leq (\vec{\tau} \wedge I\vec{\tau})(d\eta) \leq$$
$$\leq (1+\varepsilon) \sum_{j=1}^n \left( \langle d|\zeta_j|, \vec{\tau}\rangle^2 + \left\langle \frac{\eta_j}{|\zeta_j|}, \vec{\tau} \right\rangle^2 \right). \qquad (27)$$

***proof of Lemma 3.3***. The key observation is that the forms $d\eta_j = dx^j \wedge dy^j$, $\eta_j = x_1 dy^j - y_j dx^j$ and $d|\zeta_j|$ are related by

$$d\eta_j = d|\zeta_j| \wedge \frac{\eta_j}{|\zeta_j|}. \qquad (28)$$

Indeed, in metric duality with vectors, $d|\zeta_j|$ and $\frac{\eta_j}{|\zeta_j|}$ are respectively the radial unit vector in $\mathbb{C}$ and its image via $I_0$, while $d\eta_j$ is the area form. For the same reason we have, for any unit vector $v$, the equality

$$\langle d|\zeta_j|, v\rangle = \left\langle \frac{\eta_j}{|\zeta_1|}, I_0 v \right\rangle. \qquad (29)$$

By (28) and (29) we obtain

$$\langle d\eta_j, \vec{\tau} \wedge I_0 \vec{\tau}\rangle = \left\langle d|\zeta_j| \wedge \frac{\eta_j}{|\zeta_j|}, \vec{\tau} \wedge I_0 \vec{\tau} \right\rangle =$$
$$= \langle d|\zeta_j|, \vec{\tau}\rangle \left\langle \frac{\eta_j}{|\zeta_j|}, I_0 \tau \right\rangle - \left\langle \frac{\eta_j}{|\zeta_j|}, \tau \right\rangle \langle d|\zeta_j|, I_0 \vec{\tau}\rangle =$$
$$= -\langle d|\zeta_j|, \vec{\tau}\rangle \langle d|\zeta_j|, I_0 I_0 \vec{\tau}\rangle - \left\langle \frac{\eta_j}{|\zeta_j|}, \tau \right\rangle \left\langle \frac{\eta_j}{|\zeta_j|}, I_0 I_0 \tau \right\rangle =$$
$$= \langle d|\zeta_j|, \vec{\tau}\rangle^2 + \left\langle \frac{\eta_j}{|\zeta_j|}, \tau \right\rangle^2. \qquad (30)$$

This is the relation we are looking for in the case that $I$ is the standard complex structure. For a general $I$ we get a well-controlled perturbation of equality (30), namely we use (7) and by adding up (30) for $j = 1, ..., n$ we conclude the proof of the lemma.

$\square$

***Choice of good slices of*** $T$. Using coordinates $(x_1, y_1, ..., x_n, y_n, s_1, s_2)$ on $\mathbb{C}^n \times \mathbb{C}$, consider the map

$$\Lambda_\rho : \quad \mathbb{C}^n \times \mathbb{C} \to \mathbb{C}^n \times \mathbb{C} \\ (x_1, y_1, ..., x_n, y_n, s_1, s_2) \to (x_1, y_1, ..., x_n, y_n, \rho s_1, \rho s_2), \quad (31)$$

see also section 4 of [1]. This map sends (with a "vertical shrinking" with respect to the product $\mathbb{C}^n \times \mathbb{C}$) the tube $\mathbb{C}^n \times S^1_1$ to the tube $\mathbb{C}^n \times S^1_\rho$. Denote the inverse map with $(\Lambda_\rho)^{-1}$.

By slicing theory (see e.g. [8]) we have

$$\int_R^{2R} M(\langle T, |z| = r \rangle) dr \leq M(T \llcorner (B_{2R} \setminus B_R)) \leq (1+\varepsilon)\nu(0) 3\pi R^2,$$

thus at least half of the slices $\langle T, |z| = r \rangle_{r \in [R, 2R]}$ must have masses $\leq 6\pi(1+\varepsilon)\nu(0)R$. Equivalently $M\left((\Lambda_\rho)^{-1}_* \langle T, |z| = r \rangle\right) \leq 6\pi(1+\varepsilon)\nu(0)$ for at least half of the values $r \in [R, 2R]$.

Now let $\langle P, |Z| = \rho \rangle$ be a good slice and denote by $\mu$ the measure $|\text{proj}_{\vec{\tau} \wedge I\vec{\tau}}(I_0 \vec{t})| \|\langle P, |Z| = \rho \rangle\|$. We will need the following key estimate on the total variation of $\mu$.

***Estimate on*** $|\mu|$. Observe that in the tube $\mathbb{C}^n \times S^1_\rho$ we have $|I - I_0| \leq \varepsilon \rho$, since $I_0 \vec{t}$ is orthogonal to $\vec{\tau}$ and moreover $I_0 \vec{t} \cdot I_0 \vec{\tau} = \vec{t} \cdot \vec{\tau}$ (lengths and scalar product are meant as Euclidean), we have

$$|\text{proj}_{\vec{\tau} \wedge I\vec{\tau}}(I_0 \vec{t})| \leq |\vec{\tau} \cdot \vec{t}| + \varepsilon \rho.$$

Let $\nu(0)$ be the density of $T$ at 0. Denote by $\alpha_T$ the dual covector (with respect to the Euclidean metric) of $\vec{t}$. We have

$$\langle P, |Z| = \rho \rangle(\alpha_T) = \left[(\Lambda_\rho)^{-1}_* \langle P, |Z| = \rho \rangle\right](\Lambda^*_\rho \alpha_T). \quad (32)$$

By the good choice of the slice we have that $M\left((\Lambda_\rho)^{-1}_* \langle P, |Z| = \rho \rangle\right)$ is bounded by $6\pi(1+\varepsilon)\nu(0)$. The one-form $\Lambda^*_\rho \alpha_T$ on the other hand has comass $\rho$, by (31), since $\alpha_T$ is a form of unit comass only in the differentials $ds^1, ds^2$. Therefore we get

$$|\mu| \leq |\langle P, |Z| = \rho \rangle(\alpha_T)| + \varepsilon \rho \|\langle P, |Z| = \rho \rangle\| \leq \\ \leq 6\pi(1+\varepsilon)\nu(0)\rho + 6\pi(1+\varepsilon)\nu(0)\varepsilon \rho = 6\pi(1+\varepsilon)^2 \nu(0)\rho. \quad (33)$$

**Lemma 3.4.** Let $\langle P, |Z| = \rho \rangle$ be a good slice of $P \subset \mathbb{C}^n \times \mathbb{C}$ with a tube of the form $\mathbb{C}^n \times S^1_\rho$. Let $|I - I_0| \leq \frac{\varepsilon}{n}$. Denote by $\vec{\tau}$ the ($\mathcal{H}^1$-a.e. well-defined) unit vector that is the approximate tangent to $\langle P, |Z| = \rho \rangle$ and by $\vec{t}$ the unit tangent vector to the fiber $S^1_\rho$. Then we have

$$|\langle P, |Z| = \rho \rangle(\eta)| \leq (1+\varepsilon)C\nu(0)\rho \int_{\mathcal{P} \cap \{|Z|=\rho\}} \frac{(\vec{\tau} \wedge I\vec{\tau})\,(d\eta)}{|proj_{\vec{\tau} \wedge I\vec{\tau}}(I_0\vec{t})|} d\|\langle P, |Z|=\rho\rangle\|, \tag{34}$$

where $C = 3\pi$.

**proof of Lemma 3.4.** Decompose the 1-cycle $\langle P, |Z| = \rho \rangle$ as in (19).

$$\langle P, |Z| = \rho \rangle = \sum_{j=1}^{\infty} N_j \text{ with } \sum_{j=1}^{\infty} M(N_j) = M(\langle P, |Z| = r \rangle).$$

Each $N_j$ is a cycle with underlying 1-rectifiable set $\mathcal{N}_j$ that is a connected compact set. The $\mathcal{N}_j$'s are moreover disjoint. Let $w = (u_1, w_1, ..., u_n, w_n)$ be an arbitrary point in $\mathbb{R}^{2n} \equiv \mathbb{C}^n$, where we use the variables $(x_1, y_1, ..., x_n, y_n)$ in $\mathbb{R}^{2n} \equiv \mathbb{C}^n$ as already done before. The one-form $u_j dy^j - w_j dx^j$ is the exterior differential of the function $u_j y_j - w_j x_j$. Since for any $i$ the current $N_i$ is a cycle we have that (for any $j$)

$$N_i(u_j dy^j - w_j dx^j) = 0.$$

Therefore

$$N_i(\eta_j) = N_i\left(\eta_j - (u_j dy^j - w_j dx^j)\right) = N_i\left((x_j - u_j)dy^j - (y_j - w_j)dx^j\right) =$$
$$= \left[\left(\mathrm{Transl}_{(u_j, w_j)}\right)_* N_i\right](\eta_j). \tag{35}$$

In other words we can "translate horizontally" (with respect to the product $\mathbb{C}^n \times \mathbb{C}$) each $N_i$ in order to achieve the assumption that we had in Lemma 3.2, i.e. that there is a point of $\mathcal{N}_i$ that gets projected to the origin of $\mathbb{C}^n$. On the other hand the form $d\eta$ is invariant under such horizontal translations. With this in mind we have from Lemmas 3.2 and 3.3, using (35),

$$|\langle P, |Z| = \rho \rangle(\eta)| \leq \sum_{i,j} |N_i(\eta_j)| \leq (1+\varepsilon) \sum_{i=1}^{\infty} |\mu_i| \int_{\mathcal{N}_i} \frac{(\vec{\tau} \wedge I\vec{\tau})\,(d\eta)}{|\mathrm{proj}_{\vec{\tau} \wedge I\vec{\tau}}(I_0\vec{t})|} d\|N_i\|, \tag{36}$$

where we are denoting with $\mu_i$ the measure $|\mathrm{proj}_{\vec{\tau} \wedge I\vec{\tau}}(I_0\vec{t})|\,d\|N_i\|$. Since $|\mu_i| \leq |\mu|$ and the $N_i$'s add up to $\langle P, |Z| = \rho \rangle$ without overlaps we get from the latter equation

$$|\langle P, |Z| = \rho\rangle(\eta)| \leq \frac{1}{2}(1+\varepsilon)|\mu| \int_{\mathcal{P}\cap\{|Z|=\rho\}} \frac{(\vec{\tau}\wedge I\vec{\tau})(d\eta)}{|\text{proj}_{\vec{\tau}\wedge I\vec{\tau}}(I_0\vec{t})|} d\|\langle P, |Z|=\rho\rangle\|. \tag{37}$$

The proof is concluded from (37) after we recall the estimate on $\mu$ established in (33).

$\square$

*Remark* 3.3. Lemma 3.4 is, as we said in the introduction, closely related to the Poincaré inequality proved in [11]. The constant that we got on the right hand side linearly depends on $\rho$: it is the analog of the Poincaré constant on a circle of radius $\rho$ (in our case $\{0\}\times S^1_\rho$). Our proof makes this constant appear with the geometric meaning of the total variation of the measure $\mu$. The classical proof of the Poincaré inequality requires the use of the fundamental theorem of calculus in order to obtain the analogue of our inequality (22). Since we worked intrinsically (i.e. without parametrization) on the slice $\langle P, |Z|=r\rangle$, we needed, for the proof of (22), to break up the slice taking care of the connectedness. It is quite uncommon to look at such topological properties for currents and the reason why we were able to do so is that we could choose good representatives for the underlying rectifiable sets (see Remarks 3.1 and 3.2).

At this stage we are ready to prove the rate of decay in Theorem 1.1.

**proof of Theorem 1.1**. Let $R$ be arbitrary and denote by $\langle P, |Z|=\rho\rangle$ a good slice in the sense of Lemma 3.1 and of (33) with $\rho\in[R,2R]$. Keep in mind that $(P\mathbin{\llcorner}\mathcal{A}_r)(\omega_{\mathbb{C}^n})$ is a (weakly) increasing function of $r$ (Remark 2.3 showed that the action of $\omega_{\mathbb{C}^n}$ on the approximate tangents of $P$ is always non-negative). By partial integration and using respectively Lemmas 3.4 and 3.1 in the second and third line we have:

$$(P\mathbin{\llcorner}\mathcal{A}_R)(\omega_{\mathbb{C}^n}) = (P\mathbin{\llcorner}\mathcal{A}_R)(d\eta) \leq (P\mathbin{\llcorner}\mathcal{A}_\rho)(d\eta) = \langle P, |Z|=\rho\rangle(\eta) \leq$$
$$\leq C(1+\varepsilon)\nu(0)\rho \int_{\langle P,|Z|=\rho\rangle} \frac{(\vec{\tau}\wedge I\vec{\tau})(d\eta)}{|\text{proj}_{\vec{\tau}\wedge I\vec{\tau}}(I_0\vec{t})|} d\|\langle P,|Z|=\rho\rangle\| \leq$$
$$\leq C(1+\varepsilon)\nu(0)\left(P\mathbin{\llcorner}(\mathcal{A}_{2R}\setminus\mathcal{A}_R)\right)(d\eta) = C(1+\varepsilon)\nu(0)\left(P\mathbin{\llcorner}(\mathcal{A}_{2R}\setminus\mathcal{A}_R)\right)(\omega_{\mathbb{C}^n}). \tag{38}$$

This gives, by filling the hole, that for any $R>0$ it holds

$$(P\mathbin{\llcorner}\mathcal{A}_R)(\omega_{\mathbb{C}^n}) \leq \frac{C(1+\varepsilon)\nu(0)}{1+C(1+\varepsilon)\nu(0)}(P\mathbin{\llcorner}\mathcal{A}_{2R})(\omega_{\mathbb{C}^n}).$$

By a standard iteration argument we get that the following Hölder control is valid for the (non-negative) increasing function $r\to(P\mathbin{\llcorner}\mathcal{A}_r)(\omega_{\mathbb{C}^n})$:

$$(P \llcorner \mathcal{A}_r)(\omega_{\mathbb{C}^n}) \leq Kr^\gamma, \tag{39}$$

$$\text{with } K = (P \llcorner \mathcal{A}_1)(\omega_{\mathbb{C}^n}) \left(1 + \frac{1}{C(1+\varepsilon)\nu(0)}\right)$$

$$\text{and } \gamma = \log_2\left(1 + \frac{1}{C(1+\varepsilon)\nu(0)}\right),$$

with $C = 3\pi$. In particular we can see that $K$ is as small as we wish, if we have started with the current $T$ dilated enough about the origin, and

$$\frac{1}{2C(1+\varepsilon)\nu(0)} \leq \gamma \leq \frac{1}{C(1+\varepsilon)\nu(0)}.$$

Recalling Lemma 2.2 and the discussion ahead of it, the same Hölder control as in (39) is valid for the mass ratio of Theorem 1.1. $\square$

## Appendix

The decomposition result (19) was given as an easy consequence of the fact that an integral 1-cycle can be represented as a sum of Lipschitz curves. We present in this appendix a self-contained proof of (19) that generalizes to arbitrary dimensions.

**Decomposition of the slice in the style of [7]**. Recall from [7] 4.2.25 the following facts. A compactly supported integral current $C$ is called *indecomposable* if and only if there exists no integral current $R$ of the same dimension that satisfies the following two conditions simultaneously: $R \neq 0 \neq C - R$ and $M(C) + M(\partial C) = M(R) + M(\partial R) + M(C-R) + M(\partial(C-R))$. In particular remark that when $C$ is a cycle then $R$ should also be a cycle for the second condition to hold. Further we have: an arbitrary integral cycle $C$ can be written as a sum of indecomposable integral cycles $\{C_i\}_{i=1}^\infty$

$$C = \sum_{i=1}^\infty C_i \text{ with } \sum_{i=1}^\infty M(C_i) = M(C).$$

Using this result for a slice $\langle P, |Z| = r \rangle$ we get a countable family of indecomposable integral cycles $\{R_i\}_{i=1}^\infty$ such that

$$\langle P, |Z| = r \rangle = \sum_{i=1}^\infty R_i \text{ and } \sum_{i=1}^\infty M(R_i) = M(\langle P, |Z| = r \rangle). \tag{40}$$

**Decomposition of the slice taking care of the connectedness**. Consider the possible ways of writing the (compact) set $\mathcal{P} \cap \{|Z| = r\}$ as a union of the following form:

$$\mathcal{P}\cap\{|Z|=r\} = (\mathcal{P}\cap\{|Z|=r\}\cap A) \cup (\mathcal{P}\cap\{|Z|=r\}\cap B),\qquad(41)$$
$$\text{with } A \text{ and } B \text{ disjoint open sets.}$$

Unless $\mathcal{P}\cap\{|Z|=r\}$ is connected, there exists at least a decomposition of $\mathcal{P}\cap\{|Z|=r\}$ of the type (41) into two non-empty closed (and bounded) sets. Observe that, as a compact set, $\mathcal{P}\cap\{|Z|=r\}\cap A$ is a (strictly) positive distance away from the topological boundary $\partial A$.

For any choice of $A$ and $B$, the support of each $R_i$ is contained either almost completely in $A$ or almost completely in $B$. Indeed $\mathrm{supp}(R_i \llcorner A) \subset \mathcal{P}\cap\{|Z|=r\}\cap A$ and $\mathrm{supp}(R_i \llcorner B) \subset \mathcal{P}\cap\{|Z|=r\}\cap B$. Therefore the boundary of the current $R_i \llcorner A$ is supported in the intersection of $\partial A$ with the compact set $\mathcal{P}\cap\{|Z|=r\}$, and this intersection is empty. Analogously for $R_i \llcorner B$. We thus get that $R_i \llcorner A$ and $R_i \llcorner B$ are integral cycles that add up to $R_i$ (in the sense of [7] 4.2.25). By the indecomposability of $R_i$ at least one of these two cycles must be zero, which means that one of the two sets $\mathrm{supp}(R_i \llcorner A) = \mathrm{supp} R_i \cap A$, $\mathrm{supp}(R_i \llcorner B) = \mathrm{supp} R_i \cap B$ has zero $\mathcal{H}^1$-measure.

Let $\{A_\alpha\}$ and $\{B_\alpha\}$ be all the possible choices of open sets in (41) with the condition that $A_\alpha$ is always the open set containing the support of $R_1$. Then consider the sets

$$\cap_\alpha \left(\mathcal{P}\cap\{|Z|=r\}\cap A_\alpha\right),$$
$$\mathcal{P}\cap\{|Z|=r\}\cap (\cup_\alpha B_\alpha).$$

These two sets are disjoint and their union is $\mathcal{P}\cap\{|Z|=r\}$. The first set is closed bounded and connected by construction and it contains the support of $R_1$. For each $R_i$ ($i \neq 1$) we have the following dicothomy: either $R_i$ is supported completely in all of the $A_\alpha$'s or there exists an index $\overline{\alpha}$ such that the support of $R_i$ is disjoint from $A_{\overline{\alpha}}$ and lies completely in $B_{\overline{\alpha}}$. Thus each $R_i$ is supported either completely in the first set $\cap_\alpha \left(\mathcal{P}\cap\{|Z|=r\}\cap A_\alpha\right)$ or completely in the second set $\mathcal{P}\cap\{|Z|=r\}\cap (\cup_\alpha B_\alpha)$.

This means that

$$\cap_\alpha \left(\mathcal{P}\cap\{|Z|=r\}\cap A_\alpha\right) = \cup_{i\in I_1} \mathrm{supp} R_i$$

$$\mathcal{P}\cap\{|Z|=r\}\cap (\cup_\alpha B_\alpha) = \cup_{i\in \mathbb{N}\setminus I_1} \mathrm{supp} R_i$$

with $1 \in I_1$. So we have the decomposition (as cycles)

$$\langle P, |Z|=r\rangle = \sum_{i\in I_1} R_i + \sum_{i\in \mathbb{N}\setminus I_1} R_i,$$

where the first cycle on the right hand side has for support the connected compact set $\cap_\alpha \left(\mathcal{P}\cap\{|Z|=r\}\cap A_\alpha\right)$ (which contains the support of $R_1$).

By doing the same with any other $R_j$ instead of $R_1$, we decompose

$$\langle P, |Z| = r \rangle = \sum_{j=1}^{\infty} N_j \text{ with } \sum_{j=1}^{\infty} M(N_j) = M(\langle P, |Z| = r \rangle), \quad (42)$$

where each $N_j$ is an integer cycle supported in a 1-rectifiable set that is compact and connected. Write each $N_i$ as the current of integration on a 1-rectifiable set $\mathcal{N}_i$ with multiplicity $f_i \in L^1(\mathcal{N}_i, \mathbb{N})$. The decomposition (19) implies that $\cup \mathcal{N}_i$ equals $\mathcal{P} \cap \{|Z| = r\}$ a.e. and $\sum_{i=1}^{\infty} f_i = f$ a.e.

Remark however that by construction the supports of the $N_i$'s are disjoint, so we have the disjoint unions

$$\sqcup_{i=1}^{\infty} \mathrm{supp} N_i = \mathcal{P} \cap \{|Z| = r\} = \sqcup_{i=1}^{\infty} \mathcal{N}_i,$$

with the second equality to be understood in the a.e. sense. This means that for each $i$ the set $\mathcal{N}_i$ agrees $\mathcal{H}^1$-a.e. with the support of $N_i$, in particular we can choose a **good representative** $\mathcal{N}_i$ for the underlying set to the current $N_i$, namely we can choose the $\mathcal{N}_i$'s to be compact, disjoint and connected.